\documentclass[12pt]{amsart}
\usepackage{latexsym}
\usepackage[all]{xy}

\def\sw#1{{\sb{(#1)}}}

\def\<{{\langle}}
\def\>{{\rangle}}

\def\eps{\epsilon}

\def\note#1{{}}

\def\note#1{}
\def\M{{\bf M}}

\def\cC{{\mathfrak C}}

\def\hom#1#2#3#4{{{}\sb{#1}{\rm Hom}\sb{#2}(#3,#4)}}
\def\lhom#1#2#3{{{}\sb{#1}{\rm Hom}(#2,#3)}}
\def\rhom#1#2#3{{{\rm Hom}\sb{#1}(#2,#3)}}

\def\lrend#1#2#3{{{}\sb{#1}{\rm End}\sb{#2}(#3)}}

\def\beq{\begin{equation}}
\def\eeq{\end{equation}}
\def\DC{{\Delta_\cC}}
\def \eC{{\eps_\cC}}

\def\id{{I}}

\def\bdi{\begin{diagram}}
\def\edi{\end{diagram}}

\def\Label#1{\label{#1}\ifmmode\llap{[#1] }\else 
\marginpar{\smash{\hbox{\tiny [#1]}}}\fi}
\def\Label{\label}

\newtheorem{proposition}{Proposition}[section]
\newtheorem{lemma}[proposition]{Lemma}
\newtheorem{corollary}[proposition]{Corollary}
\newtheorem{theorem}[proposition]{Theorem}

\theoremstyle{definition}
\newtheorem{definition}[proposition]{Definition}
\newtheorem{example}[proposition]{Example}

\theoremstyle{remark}

\newcounter{c}

\newcommand{\etyk}[1]{\vspace{-7.4mm}$$\begin{equation}\Label{#1}
\addtocounter{c}{1}}
\renewcommand{\]}{\ifnum \value{c}=1 $$\else \end{equation}\fi}
\begin{document}

\title{Towers of corings}
\author{Tomasz Brzezi\'nski}
\address{Department of Mathematics, University of Wales Swansea,
Singleton Park, Swansea SA2 8PP, U.K.}
\email{T.Brzezinski@swansea.ac.uk}
\urladdr{http//www-maths.swan.ac.uk/staff/tb}
\date{December 2001}
\subjclass{16W30, 13B02}
\begin{abstract}
The notion of a Frobenius coring is introduced, and it is shown that any
such coring produces a tower of Frobenius corings and 
Frobenius extensions. This establishes a one-to-one correspondence 
between Frobenius corings and extensions.
\end{abstract}
\maketitle

\section{Introduction}
\noindent The notion of a coring was introduced by Sweedler in \cite{Swe:pre} as a
generalisation of coalgebras over commutative rings to the case of
non-commutative rings. Thus a coring is defined as a bimodule of a
non-commutative ring $A$, with two $A$-bimodule maps that play the role
of a coproduct and a counit. While a coalgebra can be
understood as a dualisation of an algebra, a coring is a dualisation of
a ring. Recently a number of interesting examples of
corings have been constructed, and their properties studied 
from a general point of view  in \cite{Brz:str}. This revived interest
in corings and several new results have been reported in the field (cf.\
  \cite{Abu:rat}, \cite{Brz:gro}, \cite{ElKGom:sem},
 \cite{Gom:sep}, \cite{Wis:wea},
\cite{Wis:com}). 

Corings appear naturally in the context of extensions of rings, and in
particular in the case of Frobenius extensions. 
Recall from \cite{Kas:pro} and \cite{NakTsu:fro} that a ring
 extension $A\to B$
is called a {\em Frobenius extension}
 (of the first kind) if and only if  $B$ is a
finitely generated projective right $A$-module and $B \cong 
{\rm Hom}_A(B,A)$ as $(A,B)$-bimodules. One then proves 
that a ring extension $A\to B$ is Frobenius if and only if $B$ is an
$A$-coring with a coproduct which is a $(B,B)$-bimodule map (cf.\
\cite[Proposition~4.3]{Kad:new}, \cite[Remark~2.5]{CaeIon:fro}). 

 In \cite{Mor:adj} Morita observed that an extension of rings is a 
Frobenius extension if and only if the restriction of scalars functor 
has the same right and left adjoint. Following this observation a
functor is called a 
{\em Frobenius functor} in case it
 has the same right and left adjoint (cf.\  \cite{CaeMil:Doi},  
\cite{CasGom:Fro}). This notion has been recently
extended in \cite{CaeDeG:Fro}, to cover the case of Frobenius extensions
of second kind. It has been shown in \cite{Brz:str} that a forgetful
functor from the category of right (left) comodules of an $A$-coring to the
category of right (left) $A$-modules has a right adjoint. In
this paper we focus on corings for which this forgetful functor is 
Frobenius. We show that given any such coring, one has in fact a
 family of such corings, and a family of Frobenius extensions. By this 
 means we determine which corings lead to Frobenius extensions and 
 thus provide a new way of describing  such extensions.

\subsection{Notation and conventions.} 
All rings in this paper are associative and unital, and the unit of a ring
$A$ is denoted by $1_A$, while $\M_A$  denotes the category of right
$A$-modules etc. $\cC$ is an $A$-coring, its coproduct is denoted by
$\DC$ and its counit by $\eC$. This means that $\cC$ is an 
$(A,A)$-bimodule and $\eC:\cC\to A$ and $\DC:\cC\to 
\cC\otimes_{A}\cC$ are $(A,A)$-bimodule maps such that 
$(\id_\cC\otimes_{A}\DC)\circ\DC = (\DC\otimes_{A}\id_\cC)\circ\DC$, and 
$(\id_{\cC}\otimes_{A}\eC)\circ\DC = 
(\eC\otimes_{A}\id_{\cC})\circ\DC = \id_{\cC}$. The category of 
right $\cC$-comodules is
denoted by $\M^\cC$, and a $\cC$-coaction of a right $\cC$-comodule $M$
is denoted by $\rho^M$. Recall that this means that $M$ is a right
$A$-module and $\rho^M:M\to M\otimes_A\cC$ is a right $A$-module map
such that 
$(\id_M\otimes_{A}\DC)\circ\rho^M = 
(\rho^M\otimes_{A}\id_\cC)\circ\rho^M$, and 
$(\id_{M}\otimes_{A}\eC)\circ\rho^M =\id_{M}$. We use the Sweedler notation for a coproduct,
i.e, $\DC(c) = c\sw 1\otimes_A c\sw 2$, $(\id_\cC \otimes_{A}\DC)\circ\DC 
(c) = c\sw 1\otimes_{A}c\sw 2\otimes_{A}c\sw 3$, etc.\
(summation understood) for all $c\in \cC$.
Further details on corings can be found in \cite{Brz:str}. 

\section{Frobenius corings}
\noindent The forgetful functor $F:\M^\cC\to \M_A$ has the right adjoint
$-\otimes_A\cC :\M_A\to \M^\cC$ (cf.\ \cite[Lemma~3.1]{Brz:str}). 
We study  corings
for which $-\otimes_A\cC$ is also the left adjoint of $F$, i.e., those for
which  $F$ is a Frobenius functor.
\begin{definition}
An $A$-coring $\cC$ is called a {\em Frobenius coring} provided the
forgetful functor $F:\M^\cC\to \M_A$ is a Frobenius functor.
\label{def.frob}
\end{definition}

A number of different characterisations of Frobenius corings is given in
\cite[Theorem 4.1]{Brz:str}. An important additional characterisation of
such corings has been more recently obtained in \cite[Corollary 5.6]{CaeDeG:Fro}
\begin{theorem}
$\cC$ is a Frobenius coring if and only if there exist an invariant
$e\in \cC^A = \{c\in \cC \; |\; \forall a\in A, \; ac=ca\}$ and an
$(A,A)$-bimodule map $\gamma:\cC\otimes_A\cC\to A$ such that for all
$c,c'\in \cC$,
\begin{equation}
 c\sw 1\gamma(c\sw 2\otimes_{A}c') =  \gamma(c\otimes_{A}c'\sw 
    1)c'\sw 2, \quad \gamma(c\otimes_{A}e) = \gamma(e\otimes_{A}c) = 
    \eC(c).
\label{red.fro}
\end{equation}
The pair $(\gamma, e)$ is called a {\em reduced Frobenius system} for
$\cC$.
\label{thm.cae}
\end{theorem}
As an immediate consequence of Theorem~\ref{thm.cae}  one obtains
\begin{corollary}
    Any Frobenius $A$-coring $\cC$ is a finitely generated projective 
    left and right $A$-module.
\end{corollary}
\begin{proof}
    Suppose $\cC$ is a Frobenius $A$-coring with a reduced Frobenius 
    system $(\gamma,e)$. Write $\DC(e) = \sum_{i=1}^n 
    e_i\otimes_A\bar{e}_i$. 
    Taking  $c'=e$  in  equations 
    (\ref{red.fro}) we obtain
    $c = \sum_{i=1}^{n}\gamma(c\otimes_{A}e_{i})\bar{e}_{i}$. 
    Similarly, taking $c=e$ we obtain  
    $c'= \sum_{i=1}^{n}e_{i}\gamma(\bar{e}_{i}\otimes_{A}c')$.
    Since $\gamma$ is  an $(A,A)$-bimodule map, 
    for each $i\in\{1,2,\ldots ,n\}$, the map $\xi^i : \cC\to A$, $c\mapsto 
\gamma(c\otimes_A e_i)$
is left $A$-linear while the map $\bar{\xi}^i: \cC\to A$, $c\mapsto
\gamma(\bar{e}_i\otimes_R c)$ is right $A$-linear. Hence 
$\{\xi^{i},\bar{e}_{i}\}$ is a dual basis of ${}_{A}\cC$, and 
$\{\bar{\xi}^i, e_{i}\}$ is a dual basis of $\cC_{A}$.
\end{proof}

The characterisation of Frobenius corings in Theorem~\ref{thm.cae} 
makes it easier to relate such corings to Frobenius extensions. Recall that the
statement that $A\to B$ is a Frobenius extension is equivalent to the
existence  of an $(A,A)$-bimodule map 
$E:B\to A$, and an invariant 
$\beta = \sum_{i\in I} b_i\otimes_A\bar{b}^i\in (B\otimes_AB)^B 
= \{m\in B\otimes_AB\;
|\; \forall b\in B,\; bm=mb \}$, 
such that for all $b\in B$,
\begin{equation}
\sum_{i\in I} E(bb_i)\bar{b}^i = \sum_{i\in I} b_iE(\bar{b}^ib) = b.
\label{frob.sys}
\end{equation}
$E$ is called a {\em Frobenius homomorphism}
 and  $\beta $
is known as a {\em Frobenius element}.
Frobenius element $\beta$ can
be viewed as a $(B,B)$-bimodule map $\beta: B\to B\otimes_A B$, 
$b\mapsto \beta b
= b\beta$, so that equation (\ref{frob.sys}) can be expressed as a commutative
diagram
$$
\xymatrix{
    B \ar[rr]^{\beta} \ar[d]^{\beta} \ar[drr]^{\id_{B}}& & B\otimes_AB 
    \ar[d]^{\id_{B}\otimes_{A}E} \\ 
    B\otimes_AB 
    \ar[rr]^{E\otimes_{A}\id_{B}}& & B}
$$
Now reversing the arrows in this diagram, replacing $B$ by $\cC$ and
$(B,B)$-bimodule maps by $(\cC,\cC)$-bicomodule maps, one
obtains the following commutative diagram
\begin{equation}
    \xymatrix{
    \cC \ar[rr]^{\id_{\cC}\otimes_{A}\Theta} \ar[d]^{
\Theta\otimes_{A}\id_{\cC}} 
    \ar[drr]^{\id_{\cC}} & &
    \cC\otimes_{A}\cC
    \ar[d]^{\pi} \\ 
    \cC\otimes_{A}\cC 
    \ar[rr]^{\pi}& &\cC}
\label{diag.fro}
\end{equation}
where $\pi: \cC\otimes_A\cC\to \cC$ is a $(\cC,\cC)$-bicomodule map
(with all coactions given by the coproduct $\DC$, e.g., 
$\rho^{\cC\otimes_A\cC} = \id_\cC\otimes_A \DC$) and
$\Theta:A\to \cC$ is an $(A,A)$-bimodule map. 
\begin{lemma}
Pairs $(\pi,\Theta)$ rendering commutative diagram (\ref{diag.fro})
are in one-to-one correspondence with reduced Frobenius systems $(\gamma
,e)$ for $\cC$.
\label{lemma.fro.sys}
\end{lemma}
\begin{proof}
Clearly an $(A,A)$-bimodule map $\Theta:A\to \cC$ can be identified with
$e\in \cC^A$ via $\Theta\mapsto e =\Theta(1_A)$ and $e\mapsto [a\mapsto ae=ea]$.
Under this correspondence the commutative diagrams (\ref{diag.fro})
explicitly read
\begin{equation}
\pi(c\otimes_A e) = \pi(e\otimes_A c) = c, \qquad \forall c\in \cC.
\label{pi}
\end{equation}
The
fact that $\pi:\cC\otimes_A\cC\to \cC$ is a $(\cC,\cC)$-bicomodule map
explicitly means that for all $c,c'\in \cC$,
$$
c\sw 1\otimes_A \pi(c\sw 2\otimes_{A}c') =  \DC(\pi(c\otimes_A c')) = 
\pi(c\otimes_{A}c'\sw 
    1)\otimes_A c'\sw 2
    $$
All such bicomodule maps are in one-to-one correspondence with maps
$\gamma: \cC\otimes_A \cC\to A$ satisfying the first of equations
(\ref{red.fro}) via $\pi\mapsto \gamma = \eC\circ\pi$ and $\gamma\mapsto
[\pi :c\otimes_A c' \mapsto c\sw 1\gamma(c\sw 2\otimes_A c')]$. Under
this correspondence equations (\ref{pi}) for $\pi$ translate into the 
second
and third equations in (\ref{red.fro}).
\end{proof}

In conjunction with diagram (\ref{diag.fro}), Lemma~\ref{lemma.fro.sys}
makes it clear that Frobenius corings can be viewed as dualisations of
Frobenius extensions. The pair 
$(\pi, e)$ is termed a {\em Frobenius system} for $\cC$.

\begin{proposition}
Let $A\to B$ be  a Frobenius extension with a Frobenius element $\beta$
and a Frobenius homomorphism $E$. Then $B$ is a Frobenius $A$-coring
with  coproduct $\beta$ (viewed as a $(B,B)$-bimodule map $B\to
B\otimes_A B$), counit $E$ and a Frobenius system $(\pi, 1_B)$ where
$\pi: B\otimes_A B\to B$, $b\otimes_A b'\mapsto bb'$.
\label{prop.fro.ext}
\end{proposition}
\begin{proof}
The fact that $B$ is an $A$-coring with the specified coproduct and
counit is proven in \cite[Proposition~4.3]{Kad:new},
\cite[Remark~2.5]{CaeIon:fro}, while the fact that $B$ is a Frobenius coring
 can be verified by direct calculations. We
only note that if $\beta = \sum_{i\in I} b_i\otimes_A\bar{b}^i$, 
then the fact that $\pi$ is a bicomodule morphism
means that for all $b,b'\in B$, 
$$
\sum_{i\in I}
bb_i\otimes_A\pi(\bar{b}^i\otimes_Ab') = \sum_{i\in I}
\pi(b\otimes_Ab')b_i\otimes_A\bar{b}^i = \sum_{i\in I} \pi(b\otimes_Ab_i)
\otimes_A\bar{b}^ib'.
$$
This follows immediately from the definition of $\pi$ (as a product) and from
the fact that the Frobenius element $\beta$ is $B$-central.
\end{proof}

\begin{example}
Let ${\bf End}_A$ denote the category of right
$A$-endomorphisms. Objects in ${\bf End}_A$ are pairs $(M,f_M)$, where
$M$ is a right $A$-module and $f_M$ is a right $A$-module endomorphism
of $M$. A right $A$-module map $\phi:M\to N$ is a morphism in 
${\bf End}_A$, if $\phi\circ f_M = f_N\circ \phi$. The forgetful 
functor ${\bf End}_A\to \M_A$ is a Frobenius functor. 
Its left and right adjoint $G:
\M_A\to {\bf End}_A$ is given by 
$G: M\mapsto (M,\id_M)$, $\phi \mapsto
\phi$ (all the units and counits of the adjunctions are simply identity
maps). $A$ can be viewed as a trivial $A$-coring with coproduct and
counit both identified with the identity map $I_A$. Then the category 
of comodules of $A$, 
$\M^A$ is isomorphic to  ${\bf End}_A$, and, consequently, $A$ is a
Frobenius (trivial) $A$-coring.
\end{example}

Recall from \cite{Swe:pre} that given a ring extension $B\to A$ one 
can view $\cC = A\otimes_{B}A$ as an $A$-coring with the coproduct 
$\Delta_{\cC} : A\otimes_{B} A\to  A\otimes_{B} A\otimes_{B}A$, 
$a\otimes_{B}a'\mapsto a\otimes_{B}1_A\otimes_{B}a'$ and the counit 
$\eC: a\otimes_{A}a'\mapsto aa'$. $\cC$ is known as {\em Sweedler's 
coring} associated to a ring extension $B\to A$. 
\begin{theorem}
    Let $\cC = A\otimes_{B}A$ be the Sweedler coring associated to an 
    extension $B\to A$. If $B\to A$ is a Frobenius extension then 
    $\cC$ is a Frobenius coring. Conversely, if $A$ is a faithfully 
    flat left or right $B$-module and $\cC$ is a Frobenius coring, 
    then $B\to A$ is a Frobenius extension.
\Label{thm.coFro.Swe}
\end{theorem}
\begin{proof}
    Suppose that $B\to A$ is a Frobenius extension with a Frobenius element
    $\alpha=\sum_{i\in I} a_i\otimes_B\bar{a}^i\in A\otimes_{B}A$ and a 
    Frobenius homomorphism $E: A\to B$.  Obviously $\alpha\in \cC^A$. 
    Let $e=\alpha$ and 
    $\pi = \id_A\otimes_B E \otimes_B\id_A:
A\otimes_B A \otimes_BA \cong \cC\otimes_A \cC\to \cC$. We claim that
$(\pi,e)$ is  a Frobenius
    system for $\cC$. Indeed, using the defining properties of a Frobenius
    element and a Frobenius homomorphism we have for all $a,a'\in A$,
$$
\pi(e\otimes_Aa\otimes_Ba') = \sum_{i\in I}
\pi(a_i\otimes_B\bar{a}^ia\otimes_Ba') = \sum_{i\in I}
a_iE(\bar{a}^ia)\otimes_Ba' = a\otimes_Ba'.
$$
Similarly one shows that $\pi(a\otimes_Ba'\otimes_Ae) =  a\otimes_Ba'$.
The map $\pi$ is clearly $(A,A)$-bilinear. Its $(\cC,\cC)$-bicolinearity can
easily be checked by using the fact that the image of $E$ is in $B$.

Conversely, suppose that ${}_{B}A$ or $A_{B}$ is faithfully flat, and 
that $\cC$ is a Frobenius coring with the reduced
Frobenius system $\gamma: A\otimes_BA\otimes_BA \to A$ and $e=\sum_{i\in I}
a_i\otimes_B\bar{a}^i\in \cC^A$. Using the obvious identification $\hom
AA{A\otimes_BA\otimes_BA}A \cong \lrend BBA$, view $\gamma$ as a
$(B,B)$-bimodule map $E:A\to A$. Take any $a\in A$, then 
$a=\eC(1_A\otimes_Ba) =
\gamma(1_A\otimes_Ba \otimes_A e)$, i.e., $a=\sum_{i\in I}
E(aa_i)\bar{a}^i$. Similarly one deduces that $a=\sum_{i\in I}
a_iE(\bar{a}^ia)$. Next, the first of equations (\ref{red.fro}) implies for
all $a\in A$,
$$
1_A\otimes_BE(a)  = 1_A\otimes_B \gamma(1_A\otimes_Ba\otimes_B1_A) =
\gamma(1_A\otimes_Ba\otimes_B1_A)\otimes_B1_A = E(a)\otimes_B1_A.
$$
Since $A$ is a faithfully flat left or right $B$-module, 
we conclude that $E(a)\in B$ and thus $E$ is
a Frobenius homomorphism, and $e$ is a Frobenius element.
\end{proof}

The notion of a Frobenius coring in 
Definition~\ref{def.frob} should be compared with the notions 
of right and left co-Frobenius coalgebras introduced in 
\cite{Lin:sem}. First recall that left and right duals of $\cC$ are
rings.  Explicitly, the product in the left dual 
${}^{*}\cC = \lhom A\cC 
A$ is given by $(\xi \xi')(c) = \xi(c\sw 1\xi'(c\sw 2))$ for all 
$\xi,\xi'\in {}^{*}\cC$ and all $c\in \cC$, and the 
product in the right dual $\cC^{*}= \rhom A\cC 
A$ is given by $(\xi \xi')(c) = \xi'(\xi(c\sw 1)c\sw 2)$, for all 
$\xi,\xi'\in \cC^{*}$. The units
in both cases are given by the counit $\eC$. $\cC$ is a left 
${}^{*}\cC$-module with the multiplication $\xi c = c\sw 1\xi(c\sw 
2)$, and it is a right $\cC^{*}$-module with multiplication $c\xi = 
\xi(c\sw 1)c\sw 2$. With these structures at hand, 
the notions of right or left co-Frobenius coalgebras 
can be extended to the case of corings as follows.

\begin{definition} 
    An $A$-coring $\cC$ is said to be a {\em left co-Frobenius 
    coring}
    in case there is an injective morphism $\cC\to 
    {}^{*}\cC$ of left ${}^{*}\cC$-modules. $\cC$ is a {\em right 
    co-Frobenius coring} in  case there is an injective morphism 
    $\cC\to \cC^{*}$ of right $\cC^{*}$-modules.
\Label{def.lr.coFrob}
\end{definition}

By a slight refinement of \cite[Theorem 4.1]{Brz:str} one finds 
that a Frobenius $A$-coring $\cC$ is isomorphic to ${}^{*}\cC$ as a left ${}^{*}\cC$-module 
and is isomorphic to ${\cC}^{*}$ as a right ${\cC}^{*}$-module.  
Explicitly, if $(\gamma, e)$ is a reduced Frobenius system for $\cC$, 
then the left  ${}^{*}\cC$ isomorphism is
\begin{equation}
\phi_{l}:\cC\to{}^{*}\cC, \qquad \phi_{l}: c \mapsto [c' \mapsto 
    \gamma(c'\otimes_{A}c)],
\label{iso.phi}
\end{equation}
with the inverse $\xi\mapsto   e\sw 1\xi(e\sw 2)$. The right ${\cC}^{*}$-module 
isomorphism is
$$
\phi_{r}:\cC\to{}\cC^{*}, \qquad \phi_{r}: c \mapsto [c' \mapsto 
    \gamma(c\otimes_{A}c')],
 $$
with the inverse $\xi\mapsto   \xi(e\sw 1)e\sw 2$. Thus, in particular, 
a Frobenius coring is left and right co-Frobenius.

\section{Ring structure and towers of Frobenius corings}
\noindent One of the consequences of the characterisation of 
$F$ and $-\otimes_A\cC$ as Frobenius functors in 
\cite[Theorem~4.1]{Brz:str} is that in fact $\cC$ is a Frobenius
extension of $A$. 
\begin{theorem}
    Let $\cC$ be a Frobenius $A$-coring with a  Frobenius 
    system $(\pi,e)$. Then
    
    (1) $\cC$ is a ring with product $cc' = \pi(c\otimes_{A}c')$ 
        and unit $1_{\cC}= e$.
        
     (2) $\Theta: A\to \cC$, $a\mapsto ae = ea$ is a ring map, and the
	ring extension  $\Theta: A\to \cC$ is Frobenius  with 
        a Frobenius element $\beta = \DC(e)$ and  Frobenius 
	homomorphism $E=\eC$.
\Label{cor.fro1}
\end{theorem}
\begin{proof}
Let $S = 
({}^{*}\cC)^{op}$ be the opposite ring of the left dual ring. The key 
observation here is provided by 
\cite[Theorem~4.1~(2) and (3)]{Brz:str}, which state that if $\cC$ is a 
Frobenius coring then the ring extension $\iota :A\to S $, 
$a\mapsto [c\mapsto \eC(ca)]$ is Frobenius, 
and that $\cC \cong S$ as $(A,S)$-bimodules. The isomorphism invoked 
here  is precisely the 
isomorphism $\phi_{l}$ in Equation~(\ref{iso.phi}). By 
transferring the ring structure on $S$ to $\cC$ via this isomorphism 
we obtain the ring structure on $S$ as stated in (1). This can also be 
proven directly. To this end, note that 
the alternative expressions for product 
are $cc' =  \gamma(c\otimes_{A}c'\sw 1)c'\sw 2 = c\sw 
1\gamma(c\sw 2\otimes_{A}c')$, where $\gamma = \eC\circ \pi$. 
Thus for all $c,c',c''\in\cC$ we have, 
using the left $A$-linearity of $\gamma$ and $\DC$, 
$$
(cc')c'' = (\gamma(c\otimes_{A}c'\sw 1)c'\sw 2)c'' = 
 \gamma(c\otimes_{A}c'\sw 1)\gamma(c'\sw 2 \otimes_{A}c''\sw 1)c''\sw 2.
$$
On the other hand, the first of equations (\ref{red.fro}), the 
right $A$-linearity of $\gamma$ and the left $A$-linearity of $\DC$ imply
\begin{eqnarray*}
    c(c'c'') & = & c(\gamma(c'\otimes_{A}c''\sw 1)c''\sw 2) = 
    \gamma(c\otimes_{A}\gamma(c'\otimes_{A}c''\sw 1)c''\sw 2)c''\sw 3 
    \\
    &=& \gamma(c\otimes_{A}c'\sw 1\gamma(c'\sw 2 \otimes_{A}c''\sw 1))c''\sw 2
    = \gamma(c\otimes_{A}c'\sw 1)\gamma(c'\sw 2 \otimes_{A}c''\sw 
    1)c''\sw 2. 
\end{eqnarray*}
This explicitly proves that the product in $\cC$ is associative. The 
fact that $e$ is a unit for this product follows immediately from 
equations (\ref{pi}). Note that the same ring structure is obtained 
from $(\cC^{*})^{op}$ via isomorphism $\phi_{r}$.

The ring map $\iota :A\to S$ induces  the ring map $\Theta: A\to
\cC$, which  has the form stated. Since $\iota :A\to S$ 
is a Frobenius extension, so is $\Theta:A\to \cC$. Again, this can be
verified directly. Clearly $\Theta$ is unital. It is multiplicative
since for all $a,a'\in A$ we have
$$
\Theta(a)\Theta(a') = (ae)(a'e) = \pi(ae\otimes_A a'e) = a\pi(e\otimes_A
a'e) = aa'e.
$$
By definition, the counit $\eC$ is an $(A,A)$-bimodule map. Note
further that for all $c\in \cC$,
\begin{eqnarray*}
c\DC(e) &=& ce\sw 1\otimes_{A}e\sw 2 = c\sw 1\gamma(c\sw 
2\otimes_{A}e\sw 1)\otimes_{A}e\sw 2 \\
&=& c\sw 1\otimes_{A}c\sw 2 
\gamma(c\sw 3 \otimes_{A}e) = c\sw 1\otimes_{A}c\sw 2 \eC(c\sw 3) = 
c\sw 1\otimes_{A}c\sw 2 =\DC(c),
\end{eqnarray*}
where we used equations (\ref{red.fro}) to derive the third and fourth
equalities. Similarly one shows that $\DC(c) = \DC(e)c$. 
Therefore $\DC(e)\in
(\cC\otimes_A\cC)^\cC$, and
$$
c = \eC(c\sw 1) c\sw 2 = \eC(ce\sw 1)e\sw 2, \qquad c = c\sw 1\eC(c\sw
2) = e\sw 1\eC(e\sw 2c),
$$
thus proving that $\Theta:A\to \cC$ is a Frobenius extension with a 
Frobenius element $\DC(e)$ and Frobenius homomorphism $\eC$.
\end{proof}

Note that since $A\to \cC$ is a Frobenius extension, $\cC$ has an 
$A$-coring structure provided by the Frobenius element and homomorphism, 
as explained in Proposition~\ref{prop.fro.ext}. The calculation in the
proof of Theorem~\ref{cor.fro1} reveals that $c\beta = c\DC(e) = \DC(c)$, i.e., 
this new coring structure on $\cC$ coincides with the original one.
Note also that the original $(A,A)$-bimodule structure of $\cC$ can be
viewed as being induced from $\Theta$. Indeed, for all $c\in\cC$ and $a\in A$
we have
$$
\Theta(a) c = (ae)c = \pi(ae\otimes_Ac) = a\pi(e\otimes_A c) = ac,
$$
and similarly for the right $A$-multiplication. In particular this
means that $\cC$ is an $(A,\cC)$ and $(\cC,A)$-bimodule. 

Furthermore, note that Proposition~\ref{prop.fro.ext} 
also states that $\cC$ 
is a Frobenius coring. Clearly, the corresponding
induced  Frobenius system is 
precisely the same as the original one $(\pi,e)$. All this means that a
Frobenius coring and a Frobenius extension are two different, albeit dual to
each other, descriptions of the same situation.

Now we can construct a full tower of Frobenius corings. Suppose that 
$\cC$ is a Frobenius $A$-coring with a Frobenius system $(\pi,e)$. 
Then by Theorem~\ref{cor.fro1}, $e$ viewed as a map 
$\Theta: A\to \cC$ is a Frobenius extension with the Frobenius element 
$\DC(e)$ and the Frobenius homomorphism $\eC$. Now 
Theorem~\ref{thm.coFro.Swe} implies that the Sweedler $\cC$-coring 
$\cC\otimes_{A}\cC$ is Frobenius with the Frobenius system 
$(\id_{\cC}\otimes_{A}\eC\otimes_{A}\id_{\cC}, \Delta(e))$. Then  
$\cC\otimes_{A}\cC$ is a ring with unit $\DC(e)$ and the product 
$(c\otimes_{A}c')(c''\otimes_{A}c''') = 
c\otimes_{A}\gamma(c'\otimes_{A}c'')c'''$, and the extension $\DC: \cC\to 
\cC\otimes_{A}\cC$ is Frobenius by Theorem~\ref{cor.fro1}. The 
Frobenius element explicitly reads $e\sw 1\otimes_{A}e\otimes_{A}e\sw 
2$ and the Frobenius homomorphism is $\pi$. Apply 
Theorem~\ref{thm.coFro.Swe} to deduce that Sweedler's 
$\cC\otimes_{A}\cC$-coring 
$(\cC\otimes_{A}\cC)\otimes_{\cC}(\cC\otimes_{A}\cC) \cong 
\cC\otimes_{A}\cC\otimes_{A}\cC$ is Frobenius. Iterating this 
procedure we obtain the main result of this paper
\begin{theorem}
    Let $\cC$ be a Frobenius $A$-coring, and let $\cC^{k}= 
    \cC^{\otimes_{A} k}$, $k=1,2,\ldots$, and $\cC^{0}=A$. Then there 
    is a sequence of ring maps
    $$
    \xymatrix@1{
        \cC^{0}  \ar[r]^{\Theta} & \cC^{1} \ar[r]^{\DC} & \cC^{2} 
        \ar[rr]^{\id_{\cC}\otimes_{A}\Theta\otimes_{A}\id_{\cC}}& & \cC^{3} 
        \ar[r] &\ldots}
    $$
in which for all $k=1,2,\ldots$, $\cC^{k-1}\to \cC^{k}$ is a Frobenius 
extension and $\cC^{k}$ is a Frobenius $\cC^{k-1}$-coring.
\label{thm.tower}
\end{theorem}
This tower of corings bears very close resemblance to the tower of 
rings introduced by Jones as the means for the  classification of 
subfactors of von Neumann algebras \cite{Jon:ind}. As observed by Kadison in 
\cite{Kad:Jon}, the proper framework for such towers of rings is 
provided by strongly separable extension. These are Frobenius 
extensions with some additional properties, and they lead not only to 
towers of rings but also to a family of idempotents which form a 
Temperley-Lieb algebra (the Jones idempotents). The notion of a 
strongly separable extension can be easily translated to corings, 
thus inducing the notion of a strongly coseparable coring. Suppose 
that $R$ is a commutative ring and let $A$ be an $R$-algebra. An 
$A$-coring $\cC$ is called a {\em strongly coseparable coring} 
if it is a Frobenius coring with a Frobenius system $(\pi, e)$ such 
that $\pi(\DC(e)) = ue$ and $\eC(e) = v1_{A}$ for some units 
$u,v\in R$. The ordered pair $(u:v)$ is called an {\em index} of 
$\cC$. Note that a Frobenius coring is strongly coseparable if $e$ is a
grouplike element (i.e., $\DC(e)=e\otimes_A e$, $\eC(e) = 1_A$), in which case its index is $(1_R:1_R)$. 
Following the same line of reasoning as before one proves 
that if $\cC$ is a strongly coseparable coring then the ring extension 
$A\to \cC$ is a strongly separable extension. Furthermore, in the tower 
in Theorem~\ref{thm.tower}, for all  $k=1,2,\ldots$, 
$\cC^{k-1}\to \cC^{k}$ is a strongly separable 
extension and $\cC^{k}$ is a strongly coseparable $\cC^{k-1}$-coring. 
The index of $\cC^{2k-1}$ is $(u:v)$ while the index of $\cC^{2k}$ 
is $(v:u)$.

\section*{Acknowledgements}
I would like to thank Lars Kadison, Gigel Militaru and Robert Wisbauer for 
discussions and comments, and the Engineering and Physical Sciences Research 
Council for an Advanced Fellowship.

\end{document}